\documentclass[11pt]{amsart}
\usepackage{amsmath,amscd}
\usepackage{amssymb}
\usepackage{amsthm}
\usepackage{mathrsfs}
\usepackage{hyperref}
\usepackage{calc}
             {\begin{list}{\arabic{enumi}.}{\usecounter{enumi}%
              \setlength{\labelsep}{0.5em}%
              \settowidth{\labelwidth}{(\arabic{enumi})}%
              \setlength{\leftmargin}{\labelwidth+\labelsep}}}%
             {\end{list}}

\newtheorem{Theorem}{Theorem}[section]

\newtheorem{Lemma}[Theorem]{Lemma}

\newtheorem{Proposition}[Theorem]{Proposition}

\theoremstyle{definition}

\newtheorem{Remark}[Theorem]{Remark}

\numberwithin{equation}{section}

\newcommand{\R}		{\mathbb{R}}

\title{Nonlocal inverse problem with boundary response}
\author{Tuhin Ghosh}
\address {$^{\ast\ast}$Department of Mathematics, 
	Universit\"{a}t Bielefeld.
	\newline
	\indent\: 
	E-mail:{\tt\  tghosh@math.uni-bielefeld.de}}

\begin{document}

\begin{abstract}
The problem of interest in this article is to study the (nonlocal) inverse problem of recovering a potential based on the boundary measurement associated with the fractional Schr\"{o}dinger equation.  Let $0<a<1$, and $u$ solves 
\[\begin{cases}
\left((-\Delta)^a + q\right)u = 0 \mbox{ in } \Omega\\
supp\, u\subseteq \overline{\Omega}\cup \overline{W}\\
\overline{W} \cap \overline{\Omega}=\emptyset.
\end{cases}
\]
We show that by making the exterior to boundary measurement as $\left(u|_{W}, \frac{u(x)}{d(x)^a}\big|_{\Sigma}\right)$, it is possible to  determine $q$ uniquely in $\Omega$, where $\Sigma\subseteq\partial\Omega$ be a non-empty open subset and $d(x)=d(x,\partial\Omega)$ denotes the boundary distance function.

We also discuss local characterization of the large $a$-harmonic functions in ball and its application which includes boundary unique continuation and local density result.
\end{abstract}
\maketitle
\section{Introduction \& main result }
In this paper, we address the so-called fractional Calder\'{o}n problem (see \cite{GSU}) through fractional Schr\"{o}dinger equation and study the global identifiability of the potential based on the boundary response. So far the fractional Calder\'{o}n problem remains studied (see \cite{GSU, GRSU, GLX}) based on the data measured in the exterior of the domain. In \cite{GSU}, it has been shown that one can recover the potential $q$ from the exterior measurements of the non-local Cauchy data $\left(v|_{W}, (-\Delta)^av|_{\widetilde{W}}\right)$, where $W,\widetilde{W}$ be the non-empty open subsets of $\Omega_e:=\Omega_e$, and $v$ solves
\[\begin{cases}
\left((-\Delta)^a + q\right)v = 0 \mbox{ in } \Omega\\
supp\, v\subseteq \overline{\Omega}\cup \overline{W}\\
\overline{W} \cap \overline{\Omega}=\emptyset.
\end{cases}
\]
Here we would like to address the inverse problem through introducing the exterior to boundary response map, and based on that we look for a new global uniqueness result of recovering the potential. 
 
Let $\Omega\subset\mathbb{R}^n$ $(n\geq 2$) be a smooth bounded domain, and $W\subset \Omega_e$ be another smooth domain such that $\overline{W}\cap\overline{\Omega}=\emptyset$. Let us take $q\in C^\infty_c(\Omega)$, and  $f\in C^\infty_c(W)$ which is extended by $0$ outside of $W$. Let $0<a<1$, and consider the fractional Schr\"{o}dinger equation. Let $u\in H^a(\mathbb{R}^n)$ be the solution of
\begin{equation}\label{ra2}
\begin{cases}\left((-\Delta)^a +q\right) u = 0 \quad\mbox{ in }\Omega\\
u = f\quad\mbox{ in }\Omega_e
\end{cases}
\end{equation}
Let us define
\[\mathcal{E}_a(\overline{\Omega}) = e^{+}d^a(x)C^\infty(\overline{\Omega})\]
where $e^{+}$ denotes extension by zero from $\Omega$ to $\mathbb{R}^n$, and $d$ is a $C^\infty$
function in $\overline{\Omega}$, positive in $\Omega$ and satisfying $d(x) = dist(x, \partial\Omega)$ near $\partial\Omega$. 

It follows from \cite[Section 6]{GSU} and \cite{Grubb15}, that $u\in \mathcal{E}_a(\overline{\Omega})$. That  $\frac{u(x)}{d^a(x)}\Big|_{x\in\partial\Omega}$ exists as a function in $C^\infty(\partial\Omega)$, see Subsection \ref{tpt}.

Let us define the exterior to boundary response map:
\begin{equation} \mathcal{A}_{q}: C^\infty_c(W) \to C^\infty(\partial\Omega) \mbox{ defined by } \mathcal{A}_{q}(f)=\frac{u(x)}{d^a(x)}\Big|_{\partial\Omega} \end{equation}
where $u\in\mathcal{E}_a(\overline{\Omega})$ solves \eqref{ra2}. 
 It is a well-defined map.  Later we explore its various properties of $\mathcal{A}_{q}$, see Subsection \ref{tpt}. 

In this article, we would like to see whether the map
\[ q\mapsto\mathcal{A}_q \quad\mbox{ is injective or not.}\]
This is the inverse problem we study here. 

Let us also introduce the partial boundary data problem. 
Let $\Sigma\subset\partial\Omega$ be a non-empty open set. We define the partial boundary response map
\begin{equation} \mathcal{A}^\Sigma_{q}: C^\infty_c(W) \to C^\infty(\Sigma) \mbox{ defined by } \mathcal{A}^\Sigma_{q}(f)=\frac{u(x)}{d^a(x)}\Big|_{\Sigma}. \end{equation}
Here is  our main results.
\begin{Theorem}\label{thm1}
	Let $\Sigma\subseteq\partial\Omega$ be a non-empty open subset.	Suppose for two different $q^1,q^2\in C^\infty_c(\Omega)$:
	\begin{equation}\mathcal{A}^\Sigma_{q^1}(f) = \mathcal{A}^\Sigma_{q^2}(f),\quad \forall f\in C^\infty_c(W)\end{equation}
	Then it implies $q^1=q^2$ in $\Omega$. 
\end{Theorem}
The Theorem \ref{thm1} is a global uniqueness result in
the inverse problem for the fractional Schr\"{o}dinger equation
with both partial exterior ($W\subset \Omega_e$) and boundary data ($\Sigma\subset\partial\Omega$). 

As an application, we will now solve a Robin boundary valued local inverse problem. Let $\Omega$, $W$, $q\in C^\infty_c(\Omega)$, $f\in C^\infty_c(W)$ are same as in Theorem \ref{thm1}. Let $0<a<1$, and we consider the following Robin boundary value problem for the degenerate elliptic equation  in the half space $\mathbb{R}^{n+1}_{+}=\{(y,x)\in (0,\infty)\times\mathbb{R}^n\}$: 
\begin{equation}\label{csi}\begin{cases}
\nabla\cdot \left(y^{1-2a}\nabla U(y,x) \right)=0 \quad\mbox{in }(0,\infty)\times\mathbb{R}^n\\[1mm]
\partial_y U(0,x) +q(x)U(0,x)=0 \quad\mbox{on }\Omega\\[1mm]
U(0,x)=\begin{cases}f \quad\mbox{on }W\\[1mm]
 0 \quad\mbox{on }\Omega_e\setminus \overline{W}
\end{cases}\end{cases}
\end{equation}
The above local problem is the Caffarelli-Silvestre extension \cite{CS} in the half space $\mathbb{R}^{n+1}_{+}$ of the nonlocal problem \eqref{ra2} on $\mathbb{R}^n$. It has a non-zero solution $U\in H^1(\mathbb{R}^{n+1}_{+},y^{1-2a})$. The trace $u(x)=U(0,x)\in H^a(\mathbb{R}^n)$ solves \eqref{ra2} where $U(y,x)\in H^1(\mathbb{R}^{n+1}_{+},y^{1-2a})$ solves the above problem \eqref{csi}. Note that, the Neumann derivative satisfies $\underset{y\to0}{\lim}\,y^{1-2a}\partial_yU(0,x)=(-\Delta)^au(x)$ in $\mathbb{R}^n$ in the sense of distribution.

Let $W,\widetilde{W}\subset\Omega_e$ be non-empty subset. We know from \cite{GSU,GRSU} that by measuring the partial Cauchy data $\left( U(0,x)\big|_{W}, \partial_yU(0,x)\big|_{\widetilde{W}}\right)$ we can determine $q$ uniquely in $\Omega$. See also the recent article \cite{CGRA} establishing the same result in a more general setup. 

Here we show as a direct consequence of Theorem \ref{thm1},  it is enough to measure $\left( U(0,x)\big|_{W}, \frac{U(0,x)}{d^a(x)}\big|_{\Sigma}\right)$ to determine $q$ uniquely in $\Omega$, where $d(x)=d(x,\partial\Omega)$, $x\in\Omega$ denotes the boundary distance function, and $\Sigma\subseteq\partial\Omega$ be any non-empty open subset. 

\begin{Theorem}
Let $\Omega, W, \Sigma, d(x)$ are as in Theorem \ref{thm1}. Let $U_1, U_2\in H^1(\mathbb{R}^{n+1}_{+}, y^{1-2a})$ be the solutions of the Robin boundary value problem \eqref{baq} for two different potentials $q_1,q_2\in C^\infty_c(\Omega)$ respectively, with the same partial Dirichlet data $U_1(0,x)\big|_{W}=U_2(0,x)\big|_{W}=f\in C^\infty_c(W)$. Suppose
 \[ \frac{U_1(0,x)}{d^a(x,\partial\Omega)}\Big|_{\Sigma}=\frac{U_2(0,x)}{d^a(x,\partial\Omega)}\big|_{\Sigma}, \quad\forall f\in C^\infty_c(W)\]
 then $q_1=q_2$ in $\Omega$. 
\end{Theorem}
The sole point is here we are measuring on $n$ dimensional domain $(W)$ $\times$ $(n-1)$ dimensional domain $(\Sigma)$ to recover the potential $q$ uniquely in $\Omega$, a $n$-dimensional domain. Whereas, in the previous case we were measuring on  $n$ dimensional domain $(W)$ $\times$ $n$ dimensional domain $(\widetilde{W})$ to recover the potential $q$ uniquely in $\Omega$. So our new result has the merit to reach the same conclusion, but using one-dimensional less information.

\subsection*{Literature}
These type of inverse problems are often addressed as generalized  Calder\'{o}n type inverse problems.  In the original Calder\'{o}n problem \cite{Calderon1980} the objective was to know about the internal conductivity of an object from the
static voltage and current measurements at the boundary of that object. Study of the inverse boundary value problems have a long
history, in particular, in the context of electrical impedance tomography, on seismic and medical imaging, inverse scattering problems and so on. We refer to \cite{UG} and the references therein for a survey of this topic. 

The study of fractional and nonlocal operators and its related inverse problems is a very active field in recent years. These nonlocal equations appear in modeling various problems from diffusion process \cite{AFMR}, finance \cite{SCW}, image processing \cite{GGOS}, biology \cite{MAVE} etc. See \cite{BCVE, ROSS, JBRW, CGRA} for further references.
The mathematical study of inverse problems for fractional equations (one dimensional time-space) goes back to \cite{CNYY}. For multidimensional
 space- fractional equations, in particular the fractional Calder\'{o}n problem begins with the article \cite{GSU}. Subsequent development in this particular area includes results for low regularity and stability \cite{RS1,RS4}, matrix coefficients \cite{GLX}, variable coefficient \cite{COVI1}, semilinear equations \cite{LL, LiLi2}, reconstruction from single measurement \cite{GRSU}, shape detection \cite{HL1, HL2}, local and nonlocal lower order perturbation \cite{CLR, COVI2, LiLi1, BGU}, other qualitative and quantitative analysis \cite{RAQ, RSQ, GMRA, YUH, CGRA} etc. See also the survey articles \cite{SALS, RULS} and the references therein.

This paper is organized as follows. Section 1 is the introduction. In
Section 2 we review the exterior value problem of the fractional Schr\"{o}dinger equation, and discuss the well-definedness of the exterior to boundary map $\mathcal{A}_q$. Following that, we introduce the large $a$-harmonic functions and local boundary value problem for the fractional Schr\"{o}dinger equation and discuss it in details.  In Section 3, we complete the proof of the Theorem \ref{thm1}, the solution of the inverse problem. For further interest, we add one Appendix containing local characterization of the large $a$-harmonic functions in ball and its application which includes boundary unique continuation and local density result. 
\subsection*{Acknowledgement}
The research of T.G. is supported by the Collaborative Research Center, membership no. 1283, Universit\"{a}t Bielefeld.    
\section{Preliminaries: Direct problem}
We begin with a survey of preliminaries results. 
\subsection{Fractional Laplacian and fractional Sobolev space}
Let $0<a<1$ and we recall the fractional Laplacian operator $(-\Delta)^a$ defined over the space of Schwartz class functions $\mathscr{S}(\mathbb{R}^n)$:
\begin{equation}\label{vm65}
\forall x\in\mathbb{R}^n,\quad(-\Delta)^a u(x) = \mathscr{F}^{-1} \{ \lvert{\xi}\rvert^{2a} \widehat{u}(\xi) \}, \quad u \in \mathscr{S}(\mathbb{R}^n)
\end{equation}
where $\widehat{\cdot}$ and $\mathscr{F}^{-1}$ denote the Fourier transform and its inverse respectively. 

There are many equivalent definitions of the fractional Laplacian (see \cite{KWM}). For instance, it is given by the principal value integral as  ($0 < a < 1$)
\begin{align*}
(-\Delta)^a u(x) &= C_{n,a} \ \mathrm{p.v.} \int_{\mathbb{R}^n} \frac{u(x)-u(y)}{\lvert{x-y}\rvert^{n+2a}} \,dy.
\end{align*}
where $C_{n,a}=\frac{4^{a}\Gamma(\frac{n}{2}+a)}{\pi^{n/2}\Gamma(-a)}$. 
Through out the paper, $\Gamma$ stands for the usual Gamma function. 

The fractional Laplacian extends as a bounded map:
\[
(-\Delta)^a: \ H^s(\mathbb{R}^n) \mapsto H^{s-2a}(\mathbb{R}^n)
\]
as $s \in \mathbb{R}$ and $a\in (0,1)$. Here we recall that
\[ H^s(\R^n) := \{ u\in \mathcal{S}^\prime(\mathbb{R}^n)\, |\, \langle\xi\rangle^s\widehat{u}\in L^2(\mathbb{R}^n) \} \] 
where $\langle \xi \rangle = (1 + \lvert \xi \rvert^2)^{\frac{1}{2}}$,
and $\mathcal{S}^\prime(\mathbb{R}^n)$ denotes the space of tempered distributions in $\mathbb{R}^n$.

Let $\Omega$ be an open subset of $\mathbb{R}^n$ that is either bounded with smooth boundary or equal to $\mathbb{R}^n_{+}$. $r^{+}$ stands for restriction from $\mathbb{R}^n$ to $\Omega$, $e^{+}$
stands for extension by zero from 
$\Omega$ to $\mathbb{R}^n$. $r^{-}$ and $e^{-}$ are similar for $\Omega_e := \Omega_e$. We define the following Sobolev spaces.
\begin{equation*}\begin{aligned}
\overline{H}^s(\Omega) &:= r^{+}H^s(\R^n),\\
\dot{H}^s(\overline{\Omega})&:= \{ u\in H^s(\R^n)\, |\, supp\, u\subset \overline{\Omega}\}.
\end{aligned}
\end{equation*} 
The above Sobolev spaces equipped with the norm equipped with the norm 
\[
\|u\|_{\overline{H}^s(\Omega)} = \underset{v\in H^s(\mathbb{R}^n),\, v|_{\Omega}=u}{inf}\, \|v\|_{H^s(\mathbb{R}^n)}.
\]
We recall 
\begin{align*}
L^2(\mathbb{R}^{n+1}_{+},y^{1-2a})&=\{U:\mathbb{R}^{n+1}_{+}\mapsto\mathbb{R} \mbox{ measureable, } \, y^{\frac{1-2a}{2}}U\in L^2(\mathbb{R}^{n+1}_{+})\},\\[1mm]
H^1(\mathbb{R}^{n+1}_{+},y^{1-2a})&=\{U\in L^2(\mathbb{R}^{n+1}_{+},y^{1-2a}),\,\, \nabla U\in L^2(\mathbb{R}^{n+1}_{+},y^{1-2a})\}
\end{align*}
and the trace space of $H^1(\mathbb{R}^{n+1}_{+},y^{1-2a})$ at $y=0$ is $H^a(\mathbb{R}^n)$.
\subsection{Exterior value problem}\label{vm92}
Probabilistically, the fractional Laplacian operator $(-\Delta)^a$ represents the infinitesimal generator of a symmetric $2a$-stable L\'{e}vy process in the entire space.
Here we are interested in the restriction of $(-\Delta)^a$ to a bounded domain $\Omega$.
For example, one can think of the homogeneous Dirichlet exterior value problem for the fractional Laplacian operator $\left(e.g.\, (-\Delta)^av=g\right.$ $\left.\mbox{in }\Omega\mbox{ and } v=0 \mbox{ in }\Omega_e\right)$ which represents the  infinitesimal  generator of a symmetric $2a$-stable L\'{e}vy process for which particles are killed upon leaving the domain $\Omega$. See \cite{AD}.
\subsection*{Existence, uniqueness \& stability}
Let $q\in L^{\infty}(\Omega)$, then for a given $h\in\widetilde{H}^{a}(\Omega)^{*}$ and $f\in H^{a}(\mathbb{R}^{n})$ there exist a solution (weak) $u\in H^a(\mathbb{R}^n)$ solving (see \cite{GSU})
\begin{equation}\label{vm8}
\begin{cases}
((-\Delta)^a+q)u=h & \mbox{ in }\Omega,\\
u=f & \mbox{ in }\Omega_e.
\end{cases}
\end{equation}
In addition, we have the $H^{a}$ stability estimate as follows: 
\begin{equation}\label{vm13}
\|u\|_{H^{a}(\mathbb{R}^{n})}\leq C_{n,s}\left(\|h\|_{\widetilde{H}^{a}(\Omega)^{*}}+\|f\|_{H^{a}(\mathbb{R}^{n})}\right),\end{equation}
for some constant $C_{n,a}>0$ independent of $h$ and $f$.
\subsection*{Eigenvalue condition} Let $q\in L^{\infty}(\Omega)$ is such that,  
\[((-\Delta)^a+q)w = 0 \mbox{ in }\Omega, \quad w = 0 \mbox{ in }\Omega_e
\quad\mbox{ has only $w=0$ solution}.\] 
Then the above problem \eqref{vm8} has a unique solution. 
\subsection*{Regularity}
We denote $\mathcal{E}_\mu(\overline{\Omega})$ equals to $e_\Omega d^\mu C^\infty(\overline{\Omega})$ for $\mathscr{R}(\mu)>-1$, where $d(x)$ is a smooth positive extension into $\Omega$ of $dist(x,\partial\Omega)$ near $\partial\Omega$. In general, for $\mu\in\mathbb{C}$ with $\mathscr{R}(\mu)>-1$ one has $\mathcal{E}_{\mu-k}(\overline{\Omega})= span\, D^{(k)}\mathcal{E}_\mu(\overline{\Omega})$, where $D^{(k)}$ denotes the smooth differential operators of order $k\in\mathbb{N}$. 

Let $q\in C^\infty_c(\Omega)$ and $u\in \dot{H}^a(\mathbb{R}^n)$ be the homogeneous Dirichlet problem     
\begin{equation}\label{vm1}\begin{cases}\left((-\Delta)^a + q\right)u = h \mbox{ in } \Omega\\
u = 0 \quad\mbox{in }\mathbb{R}^n\setminus \Omega.
\end{cases}
\end{equation}
Then due to the results in \cite{RSE2} and \cite{Grubb15} respectively, we have the following:  
\begin{itemize}
\item $h\in L^\infty(\Omega)\Longrightarrow u\in d^a\,C^\alpha(\overline{\Omega})$
 for some small $\alpha\in (0,1)$, with satisfying \begin{equation}\label{baq}\lvert|\frac{u(x)}{d^a(x)}\big|_{\partial\Omega}\rvert|_{C^\alpha(\partial\Omega)}\leq C_{n,a}\,\|h\|_{L^\infty(\Omega)}.\end{equation}
\item
$h\in C^\infty(\overline{\Omega})\Longleftrightarrow u\in \mathcal{E}_a(\overline{\Omega}).$
\end{itemize}
\subsection{Well-definedness  of the map $\mathcal{A}_q$}\label{tpt}
Let us recall that 
\[ \mathcal{A}_q: C^\infty_c(W) \to C^\infty(\partial\Omega),\qquad\mathcal{A}_q(f)= \frac{u(x)}{d^a(x)}\Big|_{\partial\Omega}\]
where $u\in H^a(\mathbb{R}^n)$ solves 
\begin{equation}\label{nra2}
\begin{cases}\left((-\Delta)^a +q\right) u = 0 \quad\mbox{ in }\Omega\\
u = f \quad\mbox{ in }\Omega_e
\end{cases}
\end{equation}
where $f\in C^\infty_c(W)$ and extended by $0$ outside, and $W$ be some smooth open set in $\Omega_e$, with satisfying $\overline{W}\cap\overline{\Omega}=\emptyset$.  

We find $w=(u-f)\in \dot{H}^a(\mathbb{R}^n)$ be the solution of  
\[\begin{cases}\left((-\Delta)^a +q\right) w = (-\Delta)^af\mbox{ in }\Omega\\[1mm]
w = 0 \quad\mbox{in }\mathbb{R}^n\setminus \Omega.
\end{cases}\]
Let us note that $h:=(-\Delta)^af\in C^\infty(\overline{\Omega})$, it simply follows since
\[ h(x)= C_{n,a}\int_W \frac{f(y)}{|x-y|^{n+2a}}dy, \quad x\in\overline{\Omega} \]
has no singularity as $ \overline{W}\cap\overline{\Omega} =\emptyset$. Therefore, by the above regularity result, we have $w\in \mathcal{E}_a(\overline{\Omega})$, that $\frac{w(x)}{d^a(x)}\big|_{\partial\Omega}\in C^\infty(\overline{\Omega})$. 

Since $\frac{w(x)}{d^a(x)}\big|_{\overline{\Omega}}=\frac{u(x)}{d^a(x)}\big|_{\overline{\Omega}}$, hence $u\in \mathcal{E}_a(\overline{\Omega})$ and in particular $\frac{u(x)}{d^a(x)}\big|_{\partial\Omega}\in C^\infty(\partial\Omega)$, where $u$ solves \eqref{nra2}.

This shows the map $f\mapsto \mathcal{A}_q(f)$ is well-defined from $C^\infty_c(W)$ to $C^\infty(\partial\Omega)$.
\subsection*{Boundedness of $\mathcal{A}_q$}The estimate \eqref{baq} asserts that we can think of $\mathcal{A}_q$ being a bounded map with respect to $L^\infty$-norm,  i.e.
\[ \|\mathcal{A}_q(f)\|_{L^\infty(\partial\Omega)}\leq C_{n,a}(dist(\overline{\Omega},\overline{W}))^{-a}\,\|f\|_{L^\infty(W)}\] 
where we have used the fact that since $\overline{\Omega}\cap \overline{W}=\emptyset$ and\[(dist(x,\overline{W}))^{-a}\sim \int_W\frac{1}{|x-y|^{n+2a}}dy,\quad x\in\overline{\Omega}.\]
\subsection*{Range characterization of $\mathcal{A}_q$} 
Let us consider the range of this map $\mathcal{A}_q$ to define the space 
\[ \mathcal{R}_{\partial\Omega}= \Big\{\frac{u(x)}{d^a(x)}\Big|_{\partial\Omega}\,:\, u \mbox{ solves }\eqref{nra2},\, f\in C^\infty_c(W)\Big\}.\] 
\begin{Proposition}\label{nra3}
	$\mathcal{R}_{\partial\Omega}$ is dense in $C^\infty(\partial\Omega)$. 	
\end{Proposition}
\begin{proof}
	The proof is given there in \cite[Theorem 1.3]{GSU}. 
\end{proof}
Since $\mathcal{A}_q$ is an open (not necessarily injective) linear map, and its range ($\mathcal{R}_{\partial\Omega}$) is dense, thus we can assert that its set of pre-images is also dense in the domain  of that map, i.e:
\[
Z= \{0\} \cup \{f\in C^\infty_c(W),\, \mathcal{A}_q(f)\neq 0\}\quad \mbox{ is dense in }C^\infty_c(W). 
\]
\subsection{Boundary value problem} In this subsection, we will introduce the local boundary value problem for the fractional 
Schr\"{o}dinger operator. This class of solutions are the key to solve our inverse problem.
\subsection*{Example of large $a$-harmonic function} The following functions are the example of a large $a$-harmonic function on the unit ball $B(0;1)$: (see \cite{AB})
\[
u_\sigma(x)=
\begin{cases}
c_{n,a}\frac{1}{(1-|x|^2)^\sigma} \mbox{ in }B(0;1)\\
c_{n,a+\sigma}\frac{1}{(|x|^2-1)^\sigma} \mbox{ in }\mathbb{R}^n\setminus\overline{B}(0;1) 
\end{cases},\quad \sigma\in (0,1-a)\]
solves \[(-\Delta)^au_\sigma=0\mbox{ in }B(0;1).\] 
Clearly, over the boundary $\partial B(0;1)$, $u_\sigma$ blows up to infinity.

The function
\begin{equation}\label{vm5}
u_{1-a}(x)=
\begin{cases}
c_{n,a}\frac{1}{(1-|x|^2)^{1-a}} \mbox{ in }B(0;1)\\[1mm]
0\quad\mbox{ in }\mathbb{R}^n\setminus\overline{B}(0;1) 
\end{cases}\end{equation}
solves
\[(-\Delta)^a u_{1-a} = 0 \mbox { in }B(0;1).\]
Over the boundary $\partial B(0; 1)$, we define the following limit 
\[\frac{u_{1-a}}{d(x)^{a-1}}|_{\partial B(0;1)}= \underset{|x|\to 1}{lim}\, \frac{u_{1-a}}{(1-|x|)^{a-1}}= \frac{c_{n,a}}{2^{1-a}}\neq 0\]
some non-zero constant. 
\subsection{Fractional Schr\"{o}dinger equation:}\label{vm15}
Motivated by the above example in general 
we would like to study the following problem:
\begin{equation}\label{vm3}\begin{cases}
\left((-\Delta)^a + q\right)u = 0 \mbox{ in } \Omega\\
supp\, u\subseteq \overline{\Omega}.
\end{cases}
\end{equation}
We will be more specific about the boundary condition soon. Before that in order to continue in that direction, let us introduce order-reducing operators of plus/minus type, see \cite{Grubb14, Grubb15, Grubb16}. We define
\[ \Xi^t_{\pm} :=Op(\chi^t_{\pm}) \mbox{ on }\mathbb{R}^n, \quad \chi^t_{\pm} = \left(\langle\xi^\prime\rangle \pm i\xi_n\right)^t.\] 
These symbols extend analytically in $\xi_n$ to $Im\, \xi_n \lessgtr 0$. Hence, by the Payley-Wiener theorem $\Xi^t_{\pm}$ preserve support in $\overline{\mathbb{R}}^n_{\pm}$. Then for all $s\in\mathbb{R}$ 
\[\Xi^t_{+}: \dot{H}^s(\overline{\mathbb{R}}^n_{+}) \to \dot{H}^{s-t}(\overline{\mathbb{R}}^n_{+})
,\quad r^{+}\Xi^t_{-}e^{+}: \overline{H}^s(\mathbb{R}^n_{+})\to \overline{H}^{s-t}(\mathbb{R}^n_{+}).\]
In fact, $\Xi^t_{+}$ and $r^{+}\Xi^t_{-}e^{+}$ are disjoints. The inverses are $\Xi^{-t}_{+}$ and $r^{+}\Xi^{-t}_{-}e^{+}$. 

The operators $\Xi^t_{+}$ maps $\mathcal{E}_t(\overline{\mathbb{R}}^n_{+})\cap \mathcal{E}^\prime$ to $e^{+}C^\infty(\overline{\mathbb{R}}^n_{+})$ with the property that 
\[ \begin{cases}\gamma_0\left(\Xi_{+}^t u\right) = \Gamma(t+1)\gamma_0(u/x_n^t), \qquad\quad\mbox{(Dirichlet)}\\
\gamma_1\left(\Xi_{+}^t u\right) = \Gamma(t+2)\gamma_0\left(\partial_{x_n}(u/x_n^t)\right), \quad\mbox{(Neumann)}\end{cases}
 u\in \mathcal{E}_t(\overline{\mathbb{R}}^n_{+})\cap \mathcal{E}^\prime\]
where $\gamma_0$ and $\gamma_1$ are the Dirichlet and Neumann boundary trace operator respectively. And $\Gamma$ stands as the well-known Gamma function.

Following that we define the $a(s)$-transmission spaces, which can be thought of as a generalization of $\mathcal{E}_a$ spaces. We define, as in \cite{Grubb15},
\[ H^{a(s)}(\overline{\mathbb{R}}^n_{+}) := \Xi^{-a}_{+}\overline{H}^{s-a}(\mathbb{R}^n_{+}),\quad\mbox{for }s-a>-\frac{1}{2}.\]
Here $e^{+}\overline{H}^{s-a}(\mathbb{R}^n_{+})$ generally has a jump at $x_n=0$; it is mapped by $\Xi^{-a}_{+}$ to a singularity of the type $x_n^a$. 
\begin{itemize}
\item  $H^{a(s)}(\overline{\mathbb{R}}^n_{+})\subset \dot{H}^{a-1/2}(\overline{\mathbb{R}}^n_{+})$ with continuous inclusions. \vspace{1mm}
\item 
$\dot{H}^{a-1/2}(\overline{\mathbb{R}}^n_{+}) \subseteq H^{a(s)}(\overline{\mathbb{R}}^n_{+}) \subseteq H_{loc}^t(\mathbb{R}^n_{+})$ with continuous inclusions, i.e. multiplication by any $\chi\in C^\infty_c(\mathbb{R}^n_{+})$ 
is bounded $H^{a(s)}(\overline{\mathbb{R}}^n_{+}) \mapsto \overline{H}^a(\mathbb{R}^n_{+})$. \vspace{1mm}
\item $\cap_{s-a-\frac{1}{2}>0} H^{a(s)}(\overline{\mathbb{R}}^n_{+}) =\mathcal{E}_{a}(\overline{\mathbb{R}}^n_{+})$, and $\mathcal{E}_{a}(\overline{\mathbb{R}}^n_{+})$
 is dense in $H^{a(s)}(\overline{\mathbb{R}}^n_{+})$.
\end{itemize}
Following that, we have
\[
H^{a(s)}(\overline{\mathbb{R}}^n_{+})=\begin{cases} \dot{H}^s(\overline{\mathbb{R}}^n_{+})\quad\mbox{ if }-\frac{1}{2}< s-a <\frac{1}{2}\\
\subseteq e^{+}x_n^a\overline{H}^{s-a}(\mathbb{R}^n_{+}) + \dot{H}^s(\overline{\mathbb{R}}^n_{+})\mbox{ if }s-a>\frac{1}{2},
\end{cases}\]
with $\dot{H}^s(\overline{\mathbb{R}}^n_{+})$ replaced by $\dot{H}^{s-\epsilon}(\mathbb{R}^n_{+})$ if $s-a-\frac{1}{2}\in\mathbb{N}\cup \{0\}$.  
\begin{Remark}
It has shown in \cite{Grubb15} that in the homogeneous Dirichlet problem \eqref{vm1}, we actually have 
$ h\in H^{s-2a}(\Omega) \Longleftrightarrow u\in  H^{a(s)}(\overline{\Omega})$.
\hfill\qed\end{Remark}
Now we will study the non-homogeneous boundary value problem. 
\subsection*{Dirichlet boundary value problem} 
Let $d(x)$ is a smooth positive extension into $\Omega$ of $dist(x,\partial\Omega)$ near $\partial\Omega$.
 The trace map $\gamma_{a-1,0}: u\mapsto \Gamma(a)\gamma_0(u/d^{a-1})$ from $\mathcal{E}_{a-1}(\overline{\Omega})$ 
to $C^\infty(\partial\Omega)$ extends to $\gamma_{a-1,0}: H^{(a-1)(s)}(\overline{\Omega})\to H^{s-a+\frac{1}{2}}(\partial\Omega)$ 
for $s-a+\frac{1}{2}>0$. 
It becomes bijection as 
\[ \gamma_{a-1,0}: 
 H^{(a-1)(s)}(\overline{\Omega})/H^{a(s)}(\overline{\Omega})\to H^{s-a+\frac{1}{2}}(\partial\Omega).\]
Following that, we address this non-homogeneous Dirichlet problem:
\begin{equation}\label{vm6}
\begin{cases}\left((-\Delta)^a + q\right)u = h \quad\mbox{ in } \Omega\\
\frac{u(x)}{d(x)^{a-1}} = f(x)   \quad\mbox{ on }\partial\Omega\\
u = 0 \quad\mbox{in }\Omega_e
\end{cases}
\end{equation}
where $h\in H^{s-2a}(\Omega)$ and $f\in H^{s-a+\frac{1}{2}}(\partial\Omega)$. Then we have the following result.
\begin{Proposition}[Nonhomogeneous Dirichlet problem (\cite{Grubb14}\cite{AB})]\label{vm28}
Let $q\in C^\infty_c(\Omega)$ and $s-a+\frac{1}{2}>0$. Then
\[ \{r^{+}\left((-\Delta)^a+q\right), \gamma_{a-1,0} \}:  H^{(a-1)(s)}(\overline{\Omega})\mapsto \overline{H}^{s-2a}(\Omega)\times H^{s-2a+\frac{1}{2}}(\partial\Omega)\]
is a Fredholm mapping. 
\end{Proposition}

We also have the following regularity result:
\[ \left(h, f\right) \in C^\infty(\overline{\Omega}) \times C^\infty(\partial\Omega) \Longrightarrow u\in \mathcal{E}_{a-1}(\overline{\Omega}).\]
\subsection*{Neumann boundary value problem}
One can also study the Neumann problems as follows: 
\begin{equation}\label{vm7}
\begin{cases}\left((-\Delta)^a + q\right)u = h \quad\mbox{ in } \Omega\\
\partial_\nu\left(\frac{u(x)}{d(x)^{a-1}}\right) = f(x)   \quad\mbox{ on }\partial\Omega\\
u = 0 \quad\mbox{in }\Omega_e
\end{cases}
\end{equation}
where $h\in H^{s-2a}(\Omega)$ and $f\in H^{s-a-\frac{1}{2}}(\partial\Omega)$. The function space for $f$ follows from the Taylor expansion of $\frac{u}{d^{a-1}}$ over the boundary (normalized with Gamma coefficients). 
If we denote
\[\begin{cases}
\gamma_{a-1, 0}\,u = u_0 := \Gamma(a)\gamma_0\left(\frac{u}{d^{a-1}}\right), \mbox{  (Dirichlet value)}\\
\gamma_{a-1, 1}\,u = u_1 := \Gamma(a+1)\gamma_0\left(\partial_\nu\left(\frac{u}{d^{a-1}}\right)\right), \mbox{  (Neumann value)}
\end{cases}\]  
then
\begin{equation}\label{vm30}\gamma_{a-1, 1}u = \gamma_{a, 0}\,u^\prime, \quad\mbox{ where }u^\prime= u-\frac{1}{\Gamma(a)}\,d^{a-1}u_0.\end{equation}
In particular, when $u_0=0$, then $ \gamma_{a-1, 1}\,u = \gamma_{a, 0}\,u = \Gamma(a+1)\frac{u}{d^a}$.
\begin{Proposition}[Nonhomogeneous Neumann problem \cite{Grubb14, Grubb15}]\label{vm29}
Let $q\in C^\infty_c(\Omega)$ and $s-a-\frac{1}{2}>0$. Then 
\[ \{r^{+}\left((-\Delta)^a+q\right), \gamma_{a-1,1} \}:  H^{(a-1)(s)}(\overline{\Omega})\mapsto \overline{H}^{s-2a}(\Omega)\times H^{s-2a+\frac{1}{2}}(\partial\Omega)\]
is a Fredholm mapping. 
\end{Proposition}

We also have the similar regularity result like in the previous case  
\[ \left(h, f\right) \in C^\infty(\overline{\Omega}) \times C^\infty(\partial\Omega) \Longrightarrow u\in \mathcal{E}_{a-1}(\overline{\Omega}).\]
\begin{Remark}
If  $ \frac{u(x)}{d(x)^{a-1}}|_{\partial\Omega}=0$ on $\partial\Omega$ in Proposition \ref{vm28} or in  Proposition \ref{vm29} then the solution of respective Dirichlet or the Neumann problem 
lies in the $H^{a(s)}(\overline{\Omega})$ space.  
\hfill\qed\end{Remark}
\subsection*{Integration by-parts formula}
Let $s-a-\frac{1}{2}>0$ with $0<a<1$. Then for $u \in H^{(a-1)(s)}(\overline{\Omega})$ and $v\in H^{a(s)}(\overline{\Omega})$ ($u=v=0$ in $\Omega_e$) we have (see \cite{Grubb16})
\begin{equation}\label{vm12} \int_\Omega v\,(-\Delta)^au - \int_\Omega u\,(-\Delta)^av = -\Gamma(a)\Gamma(a+1)\int_{\partial\Omega} \frac{u}{d^{a-1}}\,\frac{v}{d^a}.\end{equation}
\section{Inverse problem: Uniqueness result}
Here we now complete the proof of the Theorem \ref{thm1}.
\begin{proof}[Proof of the Theorem \ref{thm1}]
	Let $f\in C^\infty_0(W)$ and we extend it by $0$ elsewhere. Let $v^k_f=v^k\in H^a(\mathbb{R}^n)\cap \mathcal{E}_a(\overline{\Omega})$ solves	
	\begin{equation*}
	\begin{cases}\left((-\Delta)^a +q^k\right) v^k = 0 \quad\mbox{ in }\Omega\\
	v^k = f \quad\mbox{ in }\Omega_e.
	\end{cases}\quad k=1,2.
	\end{equation*}	
	Let $\Sigma\subset\partial\Omega$ be a non-empty open subset.  By our hypothesis we have
	\begin{equation}\label{ra9}
	\frac{v^1_f}{d^a}\Big|_{\Sigma} =  \frac{v^2_f}{d^a}\Big|_{\Sigma}, \quad \forall f\in H^a_0(W)
	\end{equation}
	and we would like to show that, it implies $q^1=q^2$ in $\Omega$. 
	
	Let $g\in C^\infty(\partial\Omega)$ be some non-zero function such that $supp\, g\subseteq \overline{\Sigma}$, and   
	$u^k\in\mathcal{E}_{a-1}(\overline{\Omega})$ be the non-zero solutions of the inhomogeneous Dirichlet problem:  
	\begin{equation}\label{ra1}
	\begin{cases}\left((-\Delta)^a +q^k\right) u^k = 0 \quad\mbox{ in }\Omega\\
	\frac{u^k(x)}{d(x)^{a-1}} = g \quad\mbox{ on }\partial\Omega\\
	u^k = 0 \quad\mbox{in }\Omega_e.
	\end{cases}\quad k=1,2.
	\end{equation}
	Next we use the integration by-parts formula \eqref{vm12} for $u^k \in \mathcal{E}_{a-1}(\overline{\Omega})$ and $w^k=(v^k-f)\in \mathcal{E}_a(\overline{\Omega})$ to have
	\begin{equation}\label{ra6}
	\int_\Omega w^k\,(-\Delta)^au^k - \int_\Omega u^k\,(-\Delta)^aw^k = -\Gamma(a)\Gamma(a+1)\int_{\partial\Omega} \frac{u^k}{d^{a-1}}\,\frac{w^k}{d^a},\qquad k=1,2. 
	\end{equation}
	Since 
	\[\begin{cases}\left((-\Delta)^a +q^k\right) w^k = (-\Delta)^af\mbox{ in }\Omega\\[1mm]
	\frac{w^k}{d^a}|_{\partial\Omega}=\frac{v^k}{d^a}|_{\partial\Omega},
	\end{cases}\quad k=1,2\]
	so together with the equation \eqref{ra1}, we re-write the above identity \eqref{ra6} as
	\begin{equation}\label{gov}
	\int_\Omega u^k\,(-\Delta)^af =\Gamma(a)\Gamma(a+1)\int_{\Sigma} g\,\frac{v^k}{d^a}, \quad k=1,2 
	\end{equation}
	since $supp\, g\subseteq \Sigma\subseteq\partial\Omega$.
	
	Then by using our hypothesis \eqref{ra9}, we conclude 
	\begin{equation*}
	\int_\Omega (u^1-u^2)\,(-\Delta)^af =0, \quad\forall f\in C^\infty_c(W).
	\end{equation*}
	Since as the difference $(u_1-u_2)\in H^a(\mathbb{R}^n)$, so it follows that
	\begin{equation*}
	\int_W f\,(-\Delta)^a(u^1-u^2) =0, \quad\forall f\in C^\infty_c(W)
	\end{equation*}
	and hence, 
	\[(-\Delta)^a(u^1-u^2)=0 \quad\mbox{in }W.\]
	Since $(u^1-u^2)=0$ in $W$ as well (cf. \eqref{ra1}), therefore by the unique continuation of the fractional Laplacian \cite[Theorem 1.2]{GSU}, we actually have 
	\[ u^1=u^2 \quad\mbox{in }\mathbb{R}^n.\] 
	This leads to the following identity follows from the equations in \eqref{ra1}:
	\[ (q^1-q^2)u^1 = 0 \quad\mbox{in }\Omega.  \]
	Since $u^1$ is non-zero in every open set in $\Omega$, thanks to the unique continuation of the fractional Laplacian \cite[Theorem 1.2]{GSU}, and it is smooth in the set containing the support of $q^1,q^2\in C^\infty_c(\Omega)$, therefore we conclude 
	\[ q^1\equiv q^2.\] 
	This completes our discussion of the proof. 
\end{proof}
 
\section{Appendix}
As we see our proof relies on the class of solutions of the non-local boundary valued problem, here we would like to offer some more discussion and new results about them.  
\subsection{Local characterization of the large $a$-harmonic functions in ball \& application}
It is possible to give a local characterization of the large $a$-harmonic functions in ball as the solution of local boundary value problem.  Let $u\in \mathcal{E}_{a-1}(\overline{\Omega})$ be the solution of
\[\begin{cases}(-\Delta)^au = 0 \quad\mbox{ in } \Omega\\
\frac{u(x)}{d(x)^{a-1}} = f(x)   \quad\mbox{ on }\partial\Omega\\
u = 0 \quad\mbox{in }\Omega_e.
\end{cases}\]
We will provide a local characterization of the related Poisson operator of the above boundary value problem, and then show some application including unique continuation result and density result.  
\subsection*{Dirichlet Green's kernel }  
Let us recall the Green's kernel associated with the fractional Laplacian operator in a bounded domain. Let $\Omega\subset\mathbb{R}^n$ be a Lipschitz domain, then we define the Green's kernel $G^a_\Omega(\cdot,\cdot)$ as 
\[G^a_\Omega(x, z)=c_{n,-a}\frac{1}{|x-z|^{n-2a}}-H^a_\Omega(x,z), \quad x, z\in \Omega,\quad x\neq z,\] where $H^a_\Omega\in H^a(\mathbb{R}^n)\cap C^a(\mathbb{R}^n)$  solves
\[\begin{cases}
(-\Delta)^a H^a_\Omega(x,\cdot) = 0 \mbox{ in }\Omega\\ 
H^a_\Omega(x, z)= c_{n,-a}\frac{1}{|x-z|^{n-2a}}\mbox{ in }\Omega_e.
\end{cases}
\]
Then $G^a_\Omega(x,\cdot)$ is known as the Dirichlet Green's kernel for the bounded domain $\Omega$, solving  
\[\begin{cases}
(-\Delta)^a_x G^a_\Omega(x,\cdot) = \delta_z(x) \mbox{ in }\Omega\\ 
G^a_\Omega(x, z)= 0\mbox{ in }\mathbb{R}^n\setminus \Omega.
\end{cases}
\]
Here we mention that, $G^a_\Omega(x,z)=G^a_\Omega(z,x)$ for all $z,x\in \R^n$, and as shown in \cite{RSE2}, the following limit exists and we call 
\[ \forall\,\omega\in \partial\Omega,\, x\in \Omega,\quad  D^aG_{\Omega}(x,\omega):=\underset{\Omega \ni z\to\, \omega }{lim}\,\,
\, \frac{G^a_{\Omega}(x,z)}{d^a(z)}.
\]
\subsection*{Poisson formula}
Let $n\geq 1$ and $u\in\mathcal{E}_{a-1}(\overline{\Omega})$ be the solution of the following problem
\begin{equation}\label{omegaequation}
\begin{cases}(-\Delta)^a u = 0 \quad\mbox{ in }\Omega\\
\frac{u(x)}{(d(x))^{a-1}}|_{\partial\Omega}= f\in C^\infty(\partial \Omega)\\
u = 0 \quad\mbox{in }\Omega_e.
\end{cases}
\end{equation}
By using he integration by parts formula \eqref{vm12}, the solution of \eqref{omegaequation} can be expressed as 
\begin{equation}\label{poissonformula} u(x) = \int_{\partial\Omega}  D^aG_{\Omega}(x,\omega)f(\omega)\, dS(\omega).\end{equation}

\subsection*{Case: $\Omega$ is a ball}
In particular for $\Omega=B(\theta;r)$ the ball of radius $r>0$ and centered at $\theta\in\mathbb{R}^n$, we have for $n\geq 2$ 
\begin{equation}\label{vm35} G^a_{B(\theta;r)}(x,z) =\begin{cases} 
\widetilde{c}_{n,a}\,\frac{1}{|z-x|^{n-2a}}\left(\int_0^{R_0(x,z)} \frac{t^{a-1}}{(1+t)^{n/2}}\, dt\right),\quad x,z\in B(\theta;r),\quad x\neq z,\\
0\quad\mbox{in }\mathbb{R}^n\setminus B(\theta;r)
\end{cases}\end{equation}
where 
\begin{equation}\label{vm4} R_0(x,z) = \frac{(r^2-|x-\theta|^2)(r^2-|z-\theta|^2)}{r^2|x-z|^2}\end{equation}
and $\widetilde{c}_{n,a}$ is some constant depending only on $n$ and $a$. 
\begin{Lemma}
	\begin{equation}\label{vm24}
	\forall\,\omega\in \mathbb{S}^{n-1},\,\, x\in B(0;1),\quad \underset{B(0;1) \ni z\to\, \omega }{lim}\,\,
	\, \frac{G^a_{B(0;1)}(x,z)}{(1-|z|^2)^a} =  k_n\, \frac{(1-|x|^2)^{a}}{|x-\omega|^n}
	\end{equation}
	where $\kappa_n= \frac{1}{n\alpha(n)}$,  $\alpha(n)$ is the volume of the unit ball. 
\end{Lemma}
\noindent
The proof of this Lemma also appears in the work of \cite{BK, HF}. To be self-contained we present it here.
\begin{proof}
	Let us recall the expression of $G^a_{B(0;1)}$ given in \eqref{vm35}. 
	\[ G^a_{B(0;1)}(x,z) =\begin{cases} 
	\widetilde{c}_{n,a}\,\frac{1}{|z-x|^{n-2a}}\left(\int_0^{R_0(x,z)} \frac{t^{a-1}}{(1+t)^{n/2}}\, dt\right),\quad x,z\in B(0;1),\quad x\neq z,\\
	0\quad\mbox{in }\mathbb{R}^n\setminus B(0;1)
	\end{cases}\]
	where \[ R_0(x,z) = \frac{(1-|x|^2)(1-|z|^2)}{|x-z|^2}\]
	and $\widetilde{c}_{n,a}$ is some constant depending only on $n$ and $a$.

	Note that,  for $t\in [0, R_0]$,  one simply has
	\[\frac{t^{a-1}}{(1+R_0)^{\frac{n}{2}}}\leq \frac{t^{a-1}}{(1+t)^{\frac{n}{2}}}\leq t^{a-1}.  \]
	Hence, it follows that
	\[
	\frac{1}{(1+R_0)^{\frac{n}{2}}}\frac{R_0^{a}}{a} \leq \int_0^{R_0(x,z)} \frac{t^{a-1}}{(1+t)^{n/2}}\, dt
	\leq \frac{R_0^{a}}{a}.\]
	Then by using the expression of $R_0(x;z)$, as a simple consequence of the above estimate we obtain the following limit as 
	\[ \forall\,\omega\in \mathbb{S}^{n-1},\, x\in B(0;1),\quad  \underset{B(0;1) \ni z\to\, \omega }{lim}\,\,
	\, \frac{G^a_{B(0;1)}(x,z)}{(1-|z|^2)^a}= \frac{\widetilde{c}_{n,a}}{a}\, \frac{(1-|x|^2)^{a}}{|x-\omega|^n}
	\]
	and $\frac{\widetilde{c}_{n,a}}{a}=\kappa_n= \frac{1}{n\alpha(n)}$, where $\alpha(n)$ is the volume of the unit ball. 
	This proves the lemma.  
	\hfill\end{proof}
Following that, here we present some local characterization result for a function satisfying non-local equation in ball only, and vice versa. 
\begin{Proposition}[Local characterization]\label{vm45}
	Let $n\geq 1$ and $u\in\mathcal{E}_{a-1}(\overline{B}(0;1))$ be the solution of the following problem
	\begin{equation}\label{vm44}
	\begin{cases}(-\Delta)^a u = 0 \quad\mbox{ in }B(0,1)\\
	\frac{u(x)}{(1-|x|^2)^{a-1}}|_{\partial B(0;1)}= f\in C^\infty(\partial B(0;1))\\
	u = 0 \quad\mbox{in }\mathbb{R}^n\setminus \overline{B}(0;1).
	\end{cases}
	\end{equation}
	Then $\frac{u(x)}{(1-|x|^2)^{a-1}}\in C^\infty(B(0;1))$ solves 
	\begin{equation}\label{vm48}
	\begin{cases}(-\Delta)\, \frac{u(x)}{(1-|x|^2)^{a-1}} = 0 \quad\mbox{ in }B(0,1)\\
	\frac{u(x)}{(1-|x|^2)^{a-1}}|_{\partial B(0;1)} = f\quad\mbox{on }\partial B(0;1)),
	\end{cases}
	\end{equation} 
	and vice-versa.  
\end{Proposition}
\begin{proof}

	We recall \eqref{vm24} for the expression of  $D^aG_{B(0;1)}(x,\omega)$ to write the Poisson formula \eqref{poissonformula} as
	\begin{equation}\label{vm39} u(x) = \frac{1}{n\alpha(n)}\, (1-|x|^2)^a\int_{\partial B(0;1)} \frac{f(\omega)}{|x-\omega|^n}\,dS(\omega), \quad x\in B(0;1)\end{equation}
	or,
	\[ \frac{u(x)}{(1-|x|^2)^{a-1}} =  \frac{1}{n\alpha(n)}\,(1-|x|^2)\int_{\partial B(0;1)} \frac{f(\omega)}{|x-\omega|^n}\,dS(\omega), \quad x\in B(0;1).\]
	The above expression is nothing but the Poisson integral formula of the harmonic functions, that it solves \eqref{vm48}:
	\[ (-\Delta)\, \frac{u(x)}{(1-|x|^2)^{a-1}} = 0 \mbox{ in }B(0;1),\quad  \frac{u(x)}{(1-|x|^2)^{a-1}} =  f \mbox{ on }\partial B(0;1).\]

	Conversely, let us define  $v\in C^\infty(\overline{B}(0;1))$ by the Poisson integral as 
	\[ v(x) = \frac{1}{n\alpha(n)}\, (1-|x|^2)\int_{\partial B(0;1)}\frac{f(\omega)}{|x-\omega|^n}\,dS(\omega), \quad x\in B(0;1)\]
	solving
	\[(-\Delta)v=0 \quad\mbox{in }B(0;1), \quad v= f\in C^\infty(\partial B(0;1))
	\]
	Then we define 
	\begin{align*} u(x) &=  (1-|x|^2)^{a-1}v(x),\quad x\in B(0;1)\\
	&= \frac{1}{n\alpha(n)}\, (1-|x|^2)^a\int_{\partial B(0;1)}\frac{f(\omega)}{|x-\omega|^n}\,dS(\omega)\end{align*}
	which solves \eqref{vm44} thanks to \eqref{vm39}.
	This completes the proof of the lemma. 
	\hfill\end{proof}
\begin{Remark}[Case\, $\Omega$ is half-space]
	There is a similar connection between null-solutions for $(1-\Delta
	)^a$ and $1-\Delta $ in the case where $\Omega =\mathbb{R}^n_+$ (here we
	write $x=(x',x_n)$, $x'=(x_1,\dots,x_{n-1})$). Let
	$ K_0$ be the Poisson operator with symbol $(\langle\xi
	'\rangle+i\xi _n)^{-1}$
	(i.e., $K_0 : f(x')\mapsto \mathcal{F}^{-1}_{\xi \to x}[(\langle\xi
	'\rangle+i\xi _n)^{-1}\mathcal{F}_{x'\to\xi '}f]$); it
	solves the problem
	\begin{equation}\label{1}
	\begin{cases} r^+(1-\Delta )v=0\text{ in }\mathbb{R}^n_+,\\
	\gamma _0v=f\text{ on }\mathbb{R}^{n-1},\\
	v=0\text{ in }\mathbb{R}^n_-.
	\end{cases}\end{equation}
	Moreover, as shown in detail in the Appendix of \cite{Grubb14}, the operator
	$K'_{a-1,0}=x_n^{a-1}K_0$ solves the problem
	\begin{equation}\label{2}
	\begin{cases} r^+(1-\Delta )^au=0\text{ in }\mathbb{R}^n_+,\\
	\gamma _0\big(\frac u{x_n^{a-1}}\big)=f\text{ on }\mathbb{R}^{n-1},\\
	u=0\text{ in }\mathbb{R}^n_-.
	\end{cases}\end{equation}
	(We have left out a normalizing Gamma-factor; the operator called $K_{a-1,0}$ in [Gru14] equals
	$\frac1{\Gamma
		(a)}K'_{a-1,0} = \frac1{\Gamma
		(a)}x_n^{a-1}K_0
	$. It is in later publications denoted
	$K^{a-1}_0$.)
	The operator maps  $H^{s-\frac{1}{2}}(\mathbb{R}^{n-1})$ into
	$x_n^{a-1}e^+\overline H^{s}(\mathbb{R}^n_+)\cap
	H^{(a-1)(s+a-1)}(\overline{\mathbb{R}^n_+})$ for $s>-\frac{1}{2}$, and maps $C_c^\infty ({\mathbb{R}^{n-1}})$ into $\mathcal{E}_{a-1}(\overline{\mathbb{R}^n_+})$.
	
	The remarkable fact that this Poisson-like operator for $(1-\Delta )^a$ is just
	$x_n^{a-1}$ times the Poisson operator for $1-\Delta $, allows the
	conclusion:
	\begin{center}
		\textit{ A function $u(x)$ is a solution of \eqref{2} if and
			only if $v(x)=\frac {u(x)}{x_n^{a-1}}$ solves \eqref{1}.}
	\end{center}
	This holds when $f$ is
	taken in $C_c^\infty (\mathbb{R}^{n-1})$ and accordingly
	$u\in \mathcal{E}_{a-1}(\overline{\mathbb{R}^n_+})$; and more generally when  $f\in H^{s-\frac{1}{2}}(\mathbb{R}^{n-1})$ and $u\in x_n^{a-1}e^+\overline H^{s}(\mathbb{R}^n_+)\cap H^{(a-1)(s+a-1)}(\overline{\mathbb{R}^n_+})$  with $s>-\frac{1}{2}$. 
	\hfill\qed\end{Remark}
\subsection{Application}
The above result enable us to prove the following boundary unique continuation principle in ball. 
\begin{Proposition}[Boundary UCP]\label{vm11}
	Let  $u\in \mathcal{E}_{a-1}(\overline{B}(0;1))$ be a solution of
	\begin{equation}\label{vm46}
	\begin{cases}
	(-\Delta)^a u = 0 \quad\mbox{ in }B(0,1)\\
	supp\, u\subseteq \overline{B}(0;1).
	\end{cases}
	\end{equation}
	Let $\Gamma\subset\partial B(0;1)$ be some non-empty connected open subset such that
	\[ \frac{u(x)}{(1-|x|^2)^{a-1}}|_{\Gamma} = \frac{u(x)}{(1-|x|^2)^a}|_{\Gamma}= 0.\]
	Then $u\equiv 0$.
\end{Proposition}
\begin{proof}
	Since $\frac{u(x)}{(1-|x|^2)^{a-1}}|_{\Gamma}=0$, then $\partial_{\nu}\left(\frac{u(x)}{(1-|x|^2)^{a-1}}\right)|_{\Gamma}= \frac{u(x)}{(1-|x|^2)^{a}}|_{\Gamma}$. The rest follows from the boundary unique continuation principle for the harmonic function. As we find  $(-\Delta)\, \frac{u(x)}{(1-|x|^2)^{a-1}}=0$ in $B(0;1)$ with $\frac{u(x)}{(1-|x|^2)^{a-1}}|_{\Gamma}=\partial_{\nu}\left(\frac{u(x)}{(1-|x|^2)^{a-1}}\right)|_{\Gamma}=0$. This implies $\frac{u(x)}{(1-|x|^2)^{a-1}}=0$ in $B(0;1)$; or, $u\equiv 0$ in $B(0;1)$. This completes the proof. 
	\hfill\end{proof}
\subsection*{Density result}
We will be proving the set consisting of the product of $a$-harmonic functions $\{u_1u_2\}$ in ball forms a dense set in $L^1_{loc}$.  
Let $u_1, u_2$ solving
\begin{equation}\label{vm2}
\begin{cases}(-\Delta)^a u = 0 \quad\mbox{ in }B(0,1)\\
\frac{u(x)}{(1-|x|^2)^{a-1}}|_{\partial B(0;1)}=f\in C_c^\infty(\Gamma)\\
u = 0 \quad\mbox{in }\mathbb{R}^n\setminus \overline{B}(0;1)
\end{cases}
\end{equation}
where $\Gamma\subset \partial B(0;1)$ be some non-empty open set. 
\begin{Proposition}[Density result]
Let $n\geq 2$. The set $\{u_1u_2\}$ where $u_k$ solves \eqref{vm2} is dense in $L^1_{loc}(B)$.  
\end{Proposition}
\begin{proof}
	It is enough to show that the for $h\in C_c(B(0;1))$ if
\begin{equation}\label{linear}
\int_{B(0;1)}hu_1u_2 =0 \quad\mbox{for all }u_1, u_2 \mbox{ satisfying }\eqref{vm2}
\end{equation}
then it must imply $h=0$. 	

By writing $\widetilde{h}= (1-|x|^2)^{2a-2} h$ and $v_k=\frac{u_k}{(1-|x|^2)^{a-1}}$ in the above identity \eqref{linear}, we get
\[\int_{B(0;1)} \widetilde{h}v_1v_2\, dx =0\]
where,   $v_k$ are the harmonic functions in $B(0;1)$ with $supp\, v_k|_{\partial\Omega} \subseteq \Gamma$, thanks to Proposition \ref{vm45}.

Then by using the result of the linearised Calder\'{o}n problem \cite[Theorem 1.1]{DKSU}, we conclude that $\widetilde{h}=0$, and so $h=0$. When $\Gamma=\partial B(0;1)$, the  $L^1$-density of the product of the harmonic functions was first observed by A.P. Calder\'{o}n in his seminal article \cite{Calderon1980}.   This completes our discussion of the proof.    
\end{proof}
\subsection*{Lack of injectivity}
As a final application, here we would like to roll out the following result, which says: There are non-zero $g\in L^2(\Omega)$ with $supp\, g\Subset\Omega$ and non-zero $v\in H^a(\mathbb{R}^n)$ which solves 
\begin{equation}\label{nra14}
\begin{cases}(-\Delta)^a v = g\quad\mbox{in }\Omega\\
v = 0 \quad\mbox{in }\Omega_e\\
supp\, g\Subset\Omega
\end{cases}\end{equation}
with satisfying
\begin{equation}\label{nra8}\frac{v(x)}{d(x)^{a}} = 0\mbox{ on }\partial\Omega.\end{equation}
Note that, for $a=1$, a function $v$ satisfying both \eqref{nra14} and \eqref{nra8} must be satisfying $v=0$ in $\Omega\setminus\overline{\omega}$, where $\omega=supp\, g$. This follows from the standard unique continuation result for the harmonic function, see \cite{Lax_R}. Here for $0<a<1$, we can not expect such result that $v=0$ near $\partial\Omega$, in fact that would imply $v=0$ everywhere, since it means $v=(-\Delta)^av=0$ in $\Omega\setminus\overline{\omega}$.   

Let us multiply the above equation \eqref{nra14} by $u\in \mathcal{E}_{a-1}(\overline{\Omega})$ solving 
\begin{equation}\label{nra15}
\begin{cases}(-\Delta)^a u = 0\quad\mbox{in }\Omega\\
\frac{u(x)}{d(x)^{a-1}} = f \quad\mbox{on }\partial\Omega\\
u=0 \quad\mbox{in }\Omega_e
\end{cases}\end{equation}
and then by using the integration by parts formula \eqref{vm12}, we find
\begin{equation}\label{nra16} \int_{\omega}gu =0,\quad\mbox{for all }u \mbox{ solving }\eqref{nra15}.\end{equation}

Let us take $\Omega=B(0;1)$. Since by Lemma \ref{vm45}, $u\in \mathcal{E}_{a-1}(\overline{B}(0;1))$ solving \eqref{nra15} means 
$\frac{u(x)}{(1-|x|^2)^{a-1}}\in C^\infty(\overline{B}(0;1))$ is harmonic, therefore by re-writing \eqref{nra16} as
\begin{equation}\label{nra17} \int_{\omega}g(1-|x|^2)^{a-1}\,\,\frac{u(x)}{(1-|x|^2)^{a-1}}\, dx =0,\quad\mbox{for all }\frac{u(x)}{(1-|x|^2)^{a-1}} \mbox{   harmonic in }B(0;1)\end{equation}
it implies $g(1-|x|^2)^{a-1} \in H^{\perp}(\omega)$, where $H(\omega)$ is set of all harmonic functions in $\omega$. 

Now we can choose some $0\neq h\in H^{\perp}(\omega)$ and consider   $g=h(1-|x|^2)^{1-a}\neq 0$ in \eqref{nra14}. Then we claim as the corresponding solution of \eqref{nra14}  $v\neq 0$ will be satisfying \eqref{nra8} i.e $\frac{v(x)}{d^a(x)}=0$ on $\partial B(0;1)$.  

Let us establish our claim here: In order to show \eqref{nra8} for our choice of $g$, let us multiply \eqref{nra14} by $u$ solving \eqref{nra15} where $f\in C^\infty(\partial\Omega)$ is arbitrary. Next by doing integration by parts (cf. \eqref{vm12}), we obtain 
\[ \int_{\partial B(0;1)} \frac{v(x)}{(1-|x|^2)^a}\, f \, d\sigma = \int_{\omega}g\,u \,dx\]
Since by our choice of $g$,\,  $\int_{\omega}gu =0$, (as $g(1-|x|^2)^{a-1} \in H^{\perp}(\omega)$ and $\frac{u(x)}{(1-|x|^2)^{a-1}} \in H(\omega)$), therefore $\int_{\partial B(0;1)} \frac{v(x)}{(1-|x|^2)^a} f =0$ for all $f\in C^\infty(\partial B(0;1))$, forcing   $ \frac{v(x)}{(1-|x|^2)^a}=0$ on $\partial B(0;1)$, i.e. \eqref{nra8}. 

\newpage

	\bibliographystyle{alpha}
	\bibliography{bib_frac}

\end{document}